\documentclass{SVProc}
\usepackage{makeidx}  % allows for indexgeneration
\makeindex
\usepackage{amsmath}
\usepackage{amsfonts}
\usepackage{wrapfig}
\usepackage{hyperref}
\usepackage{graphicx}
\usepackage{color}
\usepackage{url}

\begin{document}
\frontmatter          % for the preliminaries
\pagestyle{headings}  % switches on printing of running heads
\title{A Stochastic Nonlinear Dynamical System for Smoothing Noisy Eye Gaze Data}
\titlerunning{A Stochastic Nonlinear Dynamical System}  % abbreviated title (for running head)
\author{Thoa Thieu \inst{1} \and Roderick Melnik \inst{2}}
\authorrunning{Thoa Thieu and Roderick Melnik} % abbreviated author list (for running head)
\tocauthor{}
%
%\index{Ekeland, I.}
%\index{Temam, R.}
%\index{Dean, J.}
%\index{Grove, D.}
%\index{Chambers, C.}
%\index{Kim, B.}
%\index{Bertino, E.}
%
\institute{The University of Texas Rio Grande Valley, Edinburg, TX 78539, USA\\
\email{thoa.thieu@utrgv.edu}
\and
MS2Discovery Interdisciplinary Research Institute, Wilfrid Laurier University, \\75 University Ave W, Waterloo, Ontario, Canada N2L 3C5 \\ \email{rmelnik@wlu.ca}}

\maketitle              % typeset the title of the contribution

\begin{abstract}
In this study, we address the challenges associated with accurately determining gaze location on a screen, which is often compromised by noise from factors such as eye tracker limitations, calibration drift, ambient lighting changes, and eye blinks. We propose the use of an extended Kalman filter (EKF) to smooth the gaze data collected during eye-tracking experiments, and systematically explore the interaction of different system parameters. Our results demonstrate that the EKF significantly reduces noise, leading to a marked improvement in tracking accuracy. Furthermore, we show that our proposed stochastic nonlinear dynamical model aligns well with real experimental data and holds promise for applications in related fields.
\end{abstract}
\section{Introduction}
The Kalman filter is a highly effective algorithm for estimating and predicting the states of a system under uncertainty. It is widely used in various fields, including target tracking, navigation, and control systems, where accurate state estimation is crucial despite the presence of noise and other uncertainties \cite{Khodarahmi2023,Yadav2023,Wu2002neural}. The Kalman filter is the optimal linear estimator for systems described by linear models, assuming additive independent white noise in both the process and measurement dynamics. However, systems are inherently nonlinear in real-world applications such as those in biological systems, engineering, and many other fields. As a result, efforts have been made to extend the Kalman filter to handle nonlinear systems. Furthermore, to better represent real-world scenarios, it is essential to consider random fluctuations and uncertainties within the system. Inspired by these challenges, we are particularly interested in applying the Extended Kalman filter (EKF) to stochastic nonlinear dynamical systems to improve state estimation in the presence of both nonlinearity and randomness \cite{Law2015data,Toivanen2016,Yang2023outlier,Pyrhonen2023}. 

Eye gaze tracking is widely used in experimental research and user interface design. However, accurately determining a person’s gaze location on a screen in pixel coordinates is often challenging due to various sources of noise that affect the measurement process. These sources include limitations in the eye tracker’s accuracy, difficulties in maintaining calibration over time, and external factors such as changes in ambient lighting or shifts in the subject's position. Additionally, eye blinks cause intermittent interruptions in the data stream, further complicating the tracking process \cite{Niehorster2020,Belykh2023,Abaspourazad2024,Meo2024}.

In this study, we introduce an EKF to smooth gaze data collected during eye-tracking experiments. The EKF is employed within a stochastic nonlinear dynamical system that models the evolution of the gaze data as a nonlinear process with noise components. The system is specifically designed to account for the noisy nature of eye gaze measurements and to correct for intermittent data loss due to eye blinks.
Our results demonstrate that the EKF significantly improves the quality of the gaze signal by reducing high-frequency noise and compensating for lost data during blinks. Through numerical simulations, we observe that the EKF smooths the position and velocity estimates of the gaze data, resulting in a more consistent and accurate representation of the subject's gaze trajectory. Quantitative metrics (e.g., root mean square error) confirm that the EKF outperforms traditional filtering methods, such as a simple moving average (SMA), in terms of both noise reduction and accuracy.

Furthermore, the proposed stochastic nonlinear model fits well with real experimental data (e.g., \cite{Judd2009}) and provides valuable insights for applications in eye gaze tracking and other fields involving noisy time-series data (\cite{Toivanen2016,Yang2023outlier}). Usually, the accompanying plot illustrates how the model effectively captures the dynamic characteristics of the gaze signal. The EKF’s robustness to varying noise levels and its ability to handle interruptions in data make it an effective tool for eye gaze tracking and other relevant applications.

\section{Model description}
In this section, we introduce a model of stochastic dynamics in the EKF. In the context of the EKF, the system's evolution and the measurement process are both inherently stochastic. The EKF itself is built to handle nonlinear stochastic systems. It propagates the system’s state and uncertainty over time, updating the state estimate based on noisy observations.
In what follows, we discuss how this stochasticity is modeled.

\subsection{Stochastic State Transition Model} 

The state transition function $f(x_t)$, which models how the system evolves over time, incorporates process noise. In real-world systems like eye-tracking or motion tracking, uncertainty or randomness often arises from factors like sensor noise, human motion randomness, or modeling errors. This process noise is typically assumed to follow a Gaussian distribution.

The equation for state evolution at time step \( t \) (the system dynamics) can be written as:

\[
x_t = f(x_{t-1}) + w_t,
\]
where \( x_t \) represents the state at time \( t \), for example, the position and velocity of the subject. The term \( f(x_{t-1}) \) is the state transition function, describing how the state evolves from \( t-1 \) to \( t \), and \( w_t \) is process noise, which is assumed to follow a Gaussian distribution with zero mean and covariance matrix \( Q \):
\[
w_t \sim \mathcal{N}(0, Q)
\]
where \( Q \) is the process noise covariance matrix, typically a small value because the process noise usually represents small random disturbances.

In practice, the transition might represent motion dynamics, such as how the position and velocity evolve over time based on physics (e.g., constant velocity or acceleration), but with added uncertainty.

\subsection{Stochastic Measurement Model}

The measurements (or observations) of the system are typically corrupted by measurement noise. In the eye-tracking context, this noise could come from sensor inaccuracies, environmental factors, or other external disturbances. The measurement function \( h(x_t) \) relates the true state \( x_t \) to the observed data \( z_t \):

\[
z_t = h(x_t) + v_t,
\]
where \( z_t \) is the observed measurement at time \( t \) (which could be, for example, the measured position or velocity),
 \( h(x_t) \) is the measurement function, which maps the true state \( x_t \) to the observed data,
 \( v_t \) is measurement noise, which is typically Gaussian:
\[
v_t \sim \mathcal{N}(0, R)
\]
with \( R \) is the measurement noise covariance matrix.

This measurement noise adds uncertainty to the observed data and affects how the filter updates its state estimates.

\subsection{ Combining the Stochastic Dynamics in EKF}

The EKF estimates the state of the system (e.g., position and velocity) at each time step by predicting the state using the system's dynamics (state transition function) and then updating it with the measurements, considering both the process noise and measurement noise. The filter essentially tries to optimally combine these sources of uncertainty.

The two sources of stochasticity (process and measurement noise) are combined in the following way:

\subsubsection{Prediction (Propagation of State Uncertainty).}
The prediction step of the EKF accounts for process noise \( Q \). The state \( \hat{\mathbf{x}}_{t|t-1} \) is predicted based on the previous state, incorporating the system model \( f \) and process noise. The covariance matrix \( P_{t|t-1} \) is updated as well, considering how the uncertainty propagates.

\[
\hat{\mathbf{x}}_{t|t-1} = f(\hat{\mathbf{x}}_{t-1})
\]

\[
P_{t|t-1} = F_{t-1} P_{t-1} F_{t-1}^T + Q
\]
where \( F_{t-1} \) is the Jacobian of the state transition function \( f \) and \( Q \) denotes the process noise covariance, representing the uncertainty in the state transition.

\subsubsection{Update (Correcting State Estimate Based on Measurements).}
 When new measurements \( \mathbf{z}_t \) are observed, the EKF update step incorporates measurement noise  \( R \). The Kalman gain \( K_t \) is computed using the prediction covariance \( P_{t|t-1} \) and the measurement noise covariance \( R \):
 
 \[
 K_t = P_{t|t-1} H_t^T \left( H_t P_{t|t-1} H_t^T + R \right)^{-1}
 \]
 where:  \( H_t \) is the Jacobian of the observation function \( h \). The term \( R \) is the measurement noise covariance, accounting for uncertainty in the observed measurements.
 
 The updated state estimate \( \hat{\mathbf{x}}_{t|t} \) is:
 \[
 \hat{\mathbf{x}}_{t|t} = \hat{\mathbf{x}}_{t|t-1} + K_t \left( \mathbf{z}_t - h(\hat{\mathbf{x}}_{t|t-1}) \right)
 \]
 where  \( \mathbf{z}_t \) is the actual measurement at time \( t \) and \( h(\hat{\mathbf{x}}_{t|t-1}) \) is the predicted measurement based on the state estimate \( \hat{\mathbf{x}}_{t|t-1} \).
 
 Finally, the updated covariance \( P_{t|t} \) is:
 \[
 P_{t|t} = (I - K_t H_t) P_{t|t-1}
 \]

In essence, the stochastic dynamics come into play in both the prediction step (through the process noise \( Q \)) and the update step (through the measurement noise \( R \)). The filter continually balances the predicted state and the observed measurements, considering their respective uncertainties.

\subsubsection{The Role of Stochastic Dynamics in the Eye-Tracking Context.}

In general, the system involves eye-tracking data, where the subject’s position and velocity are being estimated over time. The stochastic dynamics could manifest in several ways:
\begin{itemize}
	\item[$\bullet$]  Process Noise: This could represent random fluctuations in the subject's movement, sensor inaccuracies, or unpredictable human behavior. It is modeled in the EKF by the process noise covariance matrix \( Q \), which governs how much the system state is allowed to change unpredictably.
	\item[$\bullet$] Measurement Noise: The eye-tracking sensors might have noise due to imperfections in measurement. This is reflected in the measurement model and noise covariance \( R \). The system does not directly measure the true position and velocity but noisy approximations of them.
\end{itemize}
%- Process Noise: This could represent random fluctuations in the subject's movement, sensor inaccuracies, or unpredictable human behavior. It is modeled in the EKF by the process noise covariance matrix \( Q \), which governs how much the system state is allowed to change unpredictably.
%
%- Measurement Noise: The eye-tracking sensors might have noise due to imperfections in measurement. This is reflected in the measurement model and noise covariance \( R \). The system does not directly measure the true position and velocity but noisy approximations of them.

The EKF takes both types of stochasticity into account by predicting the next state (with uncertainty) and then correcting this prediction based on noisy observations. The process noise \( Q \) allows the EKF to maintain uncertainty about the true state (e.g., position and velocity), and the measurement noise \( R \) ensures that the filter does not rely too heavily on the possibly noisy measurements.

\subsection{SMA Filtering}
%1. **Data Generation (Eye-tracking data)**
%
%The dataset (`data`) consists of eye-tracking data with time series for **position** (x and y) and **velocity** (vx and vy) for a subject.
%
%Let’s denote the eye-tracking data as:
%\[
%\mathbf{d}_t = \begin{bmatrix} \text{position}_t \\ \text{velocity}_t \end{bmatrix} = \begin{bmatrix} x_t \\ y_t \\ v_{x_t} \\ v_{y_t} \end{bmatrix}
%\]
%Where:
%- \( t \) is the time step.
%- \( x_t, y_t \) are the positions at time \( t \).
%- \( v_{x_t}, v_{y_t} \) are the velocities at time \( t \).

The SMA filter is a simple filter used to smooth time series data by averaging over a sliding window (e.g. \cite{Swari2021,Zulhakim2023}). Mathematically, for a data series \( x_t \), the moving average with window size \( w \) is:

\[
\hat{x}_t = \frac{1}{w} \sum_{i=t-w+1}^{t} x_i,
\]
where \( \hat{x}_t \) is the smoothed value at time \( t \).

%#### For your data (position and velocity), the moving average is applied as:
%\[
%\hat{\mathbf{d}}_t = \text{SMA}(\mathbf{d}_t, w)
%\]
%
%Where \( \hat{\mathbf{d}}_t \) represents the smoothed position and velocity.
%
%#### Code Implementation:
%- For each position \( x_t \) and velocity \( v_t \), apply the moving average over a window size \( w \). 
%- The moving average is applied separately to position and velocity.
%
%\[
%\mathbf{SMA}(\mathbf{d}_t) = \text{Moving Average}(x_t, y_t, v_{x_t}, v_{y_t})
%\]

\subsection{Root Mean Square Error (RMSE)}

The RMSE is not inherently a stochastic process, but it is used to evaluate the accuracy of the filter estimates (both from EKF and SMA) by comparing them to the true data. The randomness in the system will contribute to the spread or variance in the RMSE values. If the process and measurement noise are large, we would expect higher RMSE values, indicating less accurate estimates. We introduce the following RMSE for EKF and SMA:
\[
\text{RMSE} = \sqrt{\frac{1}{N} \sum_{t=1}^{N} (\mathbf{d}_t - \hat{\mathbf{d}}_t)^2},
\]
where \( \mathbf{d}_t \) is the true data and \( \hat{\mathbf{d}}_t \) is either the EKF or SMA estimate. The stochastic nature of the system introduces variability in the estimates, which is reflected in the RMSE.

%4. Summary of the Stochastic Dynamics in the Code
%
%To sum it up mathematically:
%
%- State Transition Model (with process noise):
%\[
%x_t = f(x_{t-1}) + w_t, \quad w_t \sim \mathcal{N}(0, Q)
%\]
%
%- Measurement Model (with measurement noise):
%\[
%z_t = h(x_t) + v_t, \quad v_t \sim \mathcal{N}(0, R)
%\]
%
%- EKF Prediction and Update:
%1. Predict state and covariance using the transition model and process noise \( Q \).
%2. Update the state estimate based on the measurement, using the measurement model and measurement noise \( R \).

The stochastic dynamics arise from the assumption that both the system evolution and the measurements are noisy, and the EKF filters these out to estimate the true state (position and velocity) over time. In the next section, we will analyze the proposed nonlinear stochastic dynamics and discuss some numerical examples.
\section{Numerical results}
In this section, we present two examples: one using synthetic data and the other with real eye-gaze data, to compare the performance of the EKF and SMA filters. The simulations in this section were conducted in Python. The key numerical results demonstrate that the EKF effectively reduces noise, resulting in a significant improvement in tracking accuracy in both examples.
\subsection{Example 1}
In this simulation, we will compare the theoretical estimates of position and velocity with the filtered estimates produced by the EKF and the SMA filter. The goal is to assess the performance of both filtering techniques in estimating the state variables and to understand how well the filtered estimates align with the theoretical values. We create an example of our data as follows. 
Let \( t = 0, 1, 2, \dots, n-1 \) represent the time steps.

\textbf{True Position \( p(t) \)}: The true position at each time step is given by:

\[
p(t) = \sin(t) + \epsilon_{\text{pos}}(t)
\]

where:
\begin{itemize}
	\item[$\bullet$] \( \sin(t) \) is the ideal true position based on a sine wave (representing periodic or oscillatory motion),
	\item[$\bullet$] \( \epsilon_{\text{pos}}(t) \sim \mathcal{N}(0, \sigma_{\text{pos}}^2) \) is Gaussian noise with mean 0 and standard deviation \( \sigma_{\text{pos}} \), added to simulate measurement noise. In this case, \( \sigma_{\text{pos}} = 0.1 \).
\end{itemize}

\textbf{True Velocity \( v(t) \)}: The true velocity at each time step is given by:

\[
v(t) = \cos(t) + \epsilon_{\text{vel}}(t)
\]

where:
\begin{itemize}
	\item[$\bullet$] \( \cos(t) \) is the ideal true velocity based on a cosine wave (which is the derivative of the sine wave with respect to time),
	\item[$\bullet$] \( \epsilon_{\text{vel}}(t) \sim \mathcal{N}(0, \sigma_{\text{vel}}^2) \) is Gaussian noise with mean 0 and standard deviation \( \sigma_{\text{vel}} \), added to simulate measurement noise. In this case, \( \sigma_{\text{vel}} = 0.1 \).
\end{itemize}

\textbf{Combined Data Matrix}: The data at each time step is represented as a 2D matrix of position and velocity:

\[
\text{data}[t] = \begin{bmatrix} p(t) \\ v(t) \end{bmatrix}
\]

The full data set for \( n \) time steps can be represented as:

\[
\text{data} = \begin{bmatrix}
	p(0) & v(0) \\
	p(1) & v(1) \\
	\vdots & \vdots \\
	p(n-1) & v(n-1)
\end{bmatrix}.
\]

\begin{figure}[h!]
	\centering
	\begin{tabular}{ll}
		\includegraphics[width=0.5\textwidth]{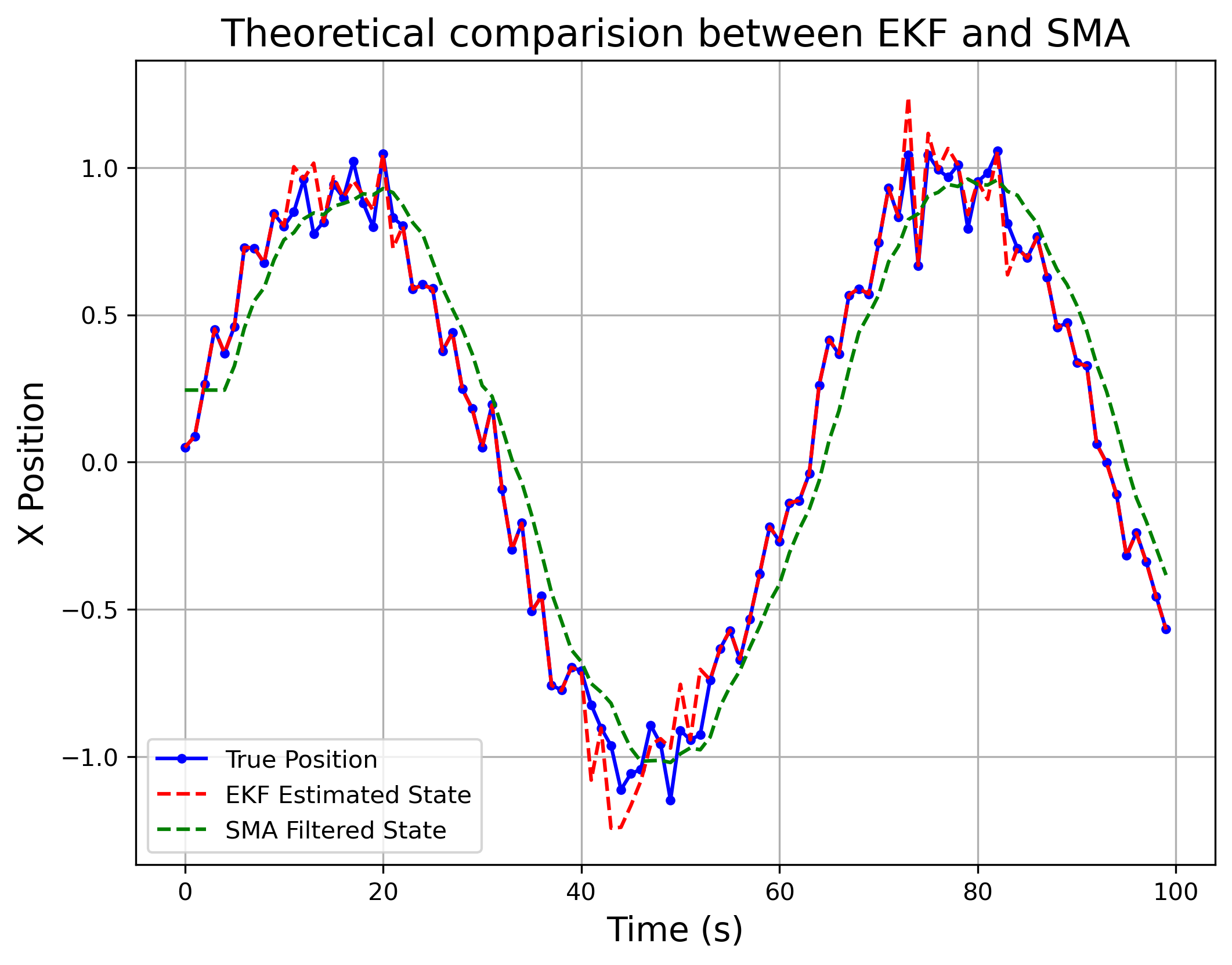}  &
		\includegraphics[width=0.5\textwidth]{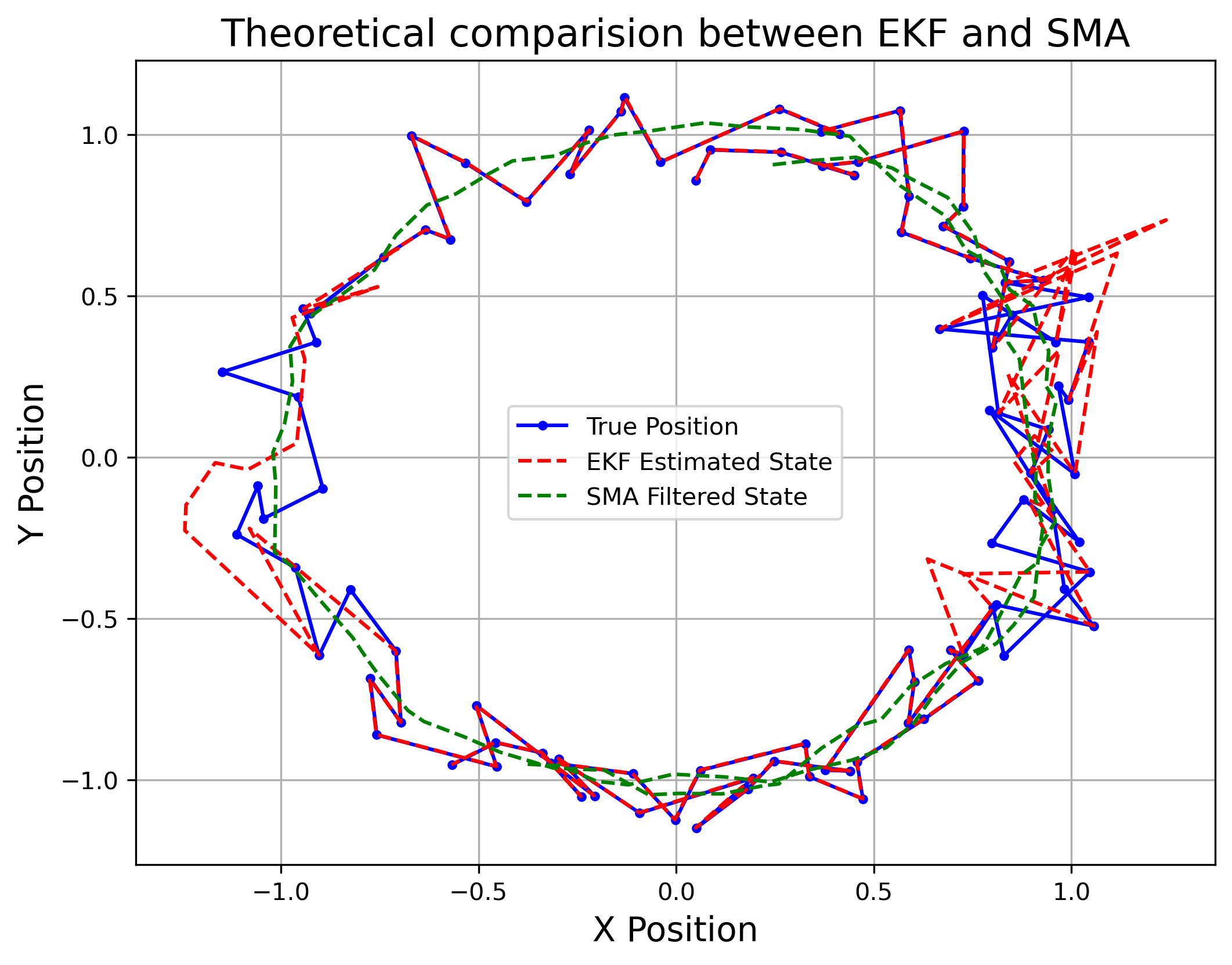}  
	\end{tabular}
	\caption{\small Comparison of EKF and SMA filters in estimating position and velocity: theoretical versus filtered estimates.
	}\label{fig:1}	
\end{figure}

Here, the first column represents the position at each time step, while the second column indicates the velocity at each time step. The main results of this example are shown in Fig. \ref{fig:1}. 

In particular, for the synthetic data, the RMSE values for the EKF and SMA are as follows:  
\begin{itemize}
	\item[$\bullet$]  RMSE (EKF) = 0.07 for position and 0.09 for velocity.  
	\item[$\bullet$] RMSE (SMA) = 0.16 for position and 0.16 for velocity.  
\end{itemize}

These results demonstrate that the EKF consistently outperforms the SMA in both position and velocity estimation. The EKF’s lower RMSE values indicate that it is able to more accurately track the true position and velocity of the system, even in the presence of noise. The SMA, on the other hand, smooths the data by averaging over a window, which introduces a lag and reduces its ability to follow rapid changes in the data. Moreover, the EKF’s lower RMSE values reflect its ability to account for the dynamic changes in the system's state through its recursive update process, which adapts to new measurements as they arrive. In contrast, the SMA is less sensitive to changes and tends to smooth out fine details, leading to a higher RMSE, particularly for rapidly changing data.
\subsection{Example 2}
\begin{figure}[h!]
	\centering
	\begin{tabular}{ll}
		\includegraphics[width=0.5\textwidth]{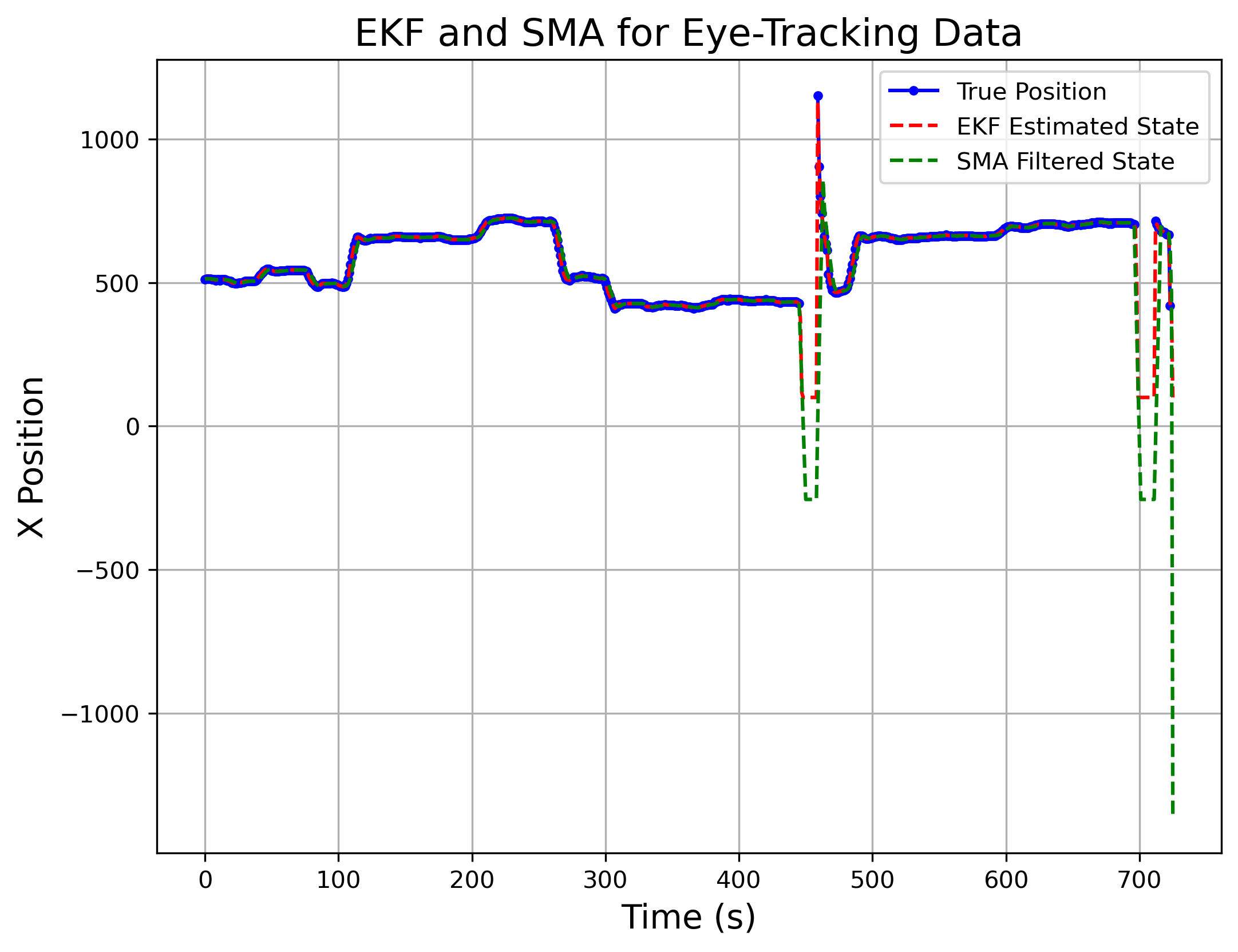}  &
		\includegraphics[width=0.5\textwidth]{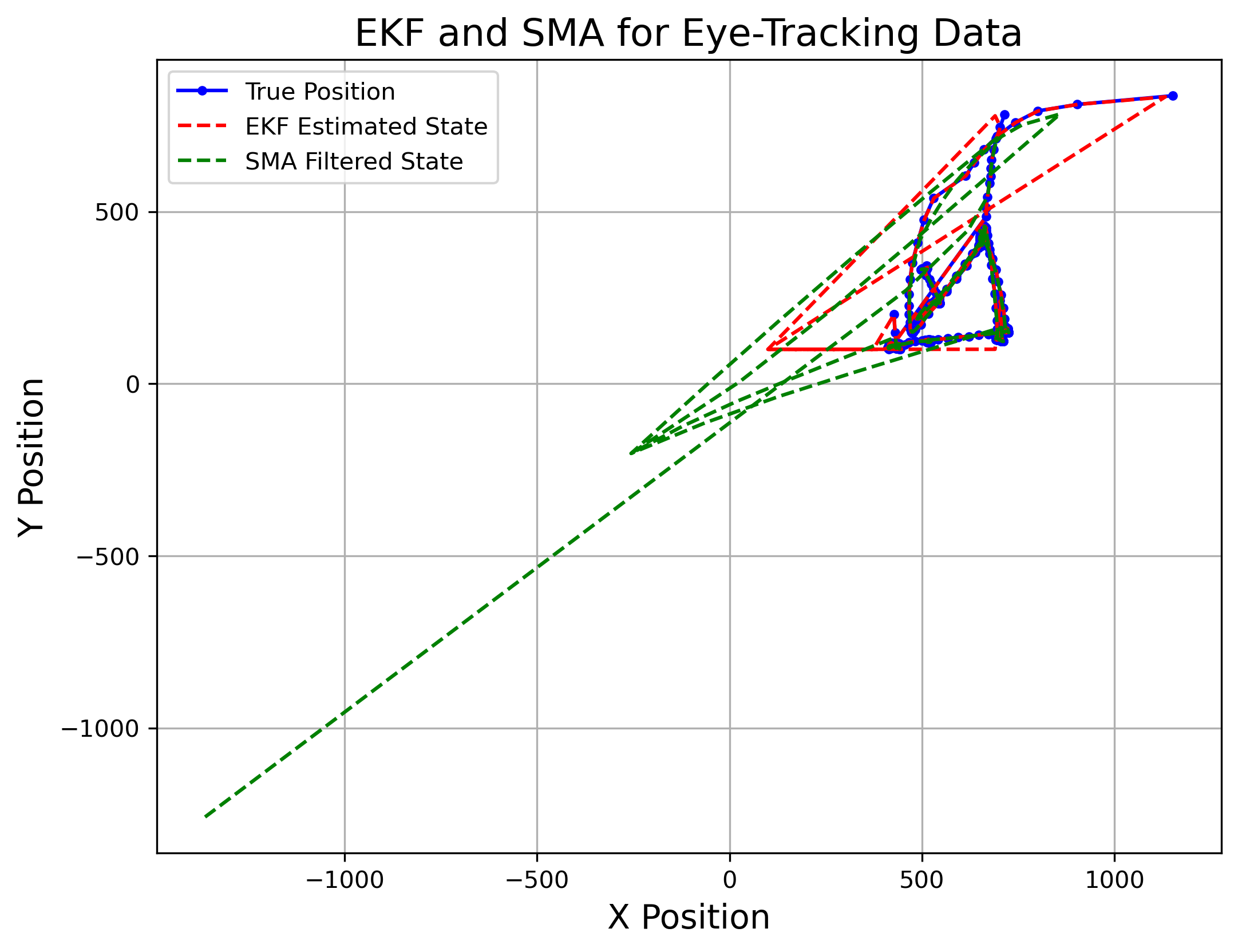}  
	\end{tabular}
	\caption{\small Comparison of EKF and SMA filters in estimating position and velocity: real-eye gaze tracking data and filtered estimates. Here, we use the data from \cite{Judd2009}.
		%		TThe circuit diagram (left) represents the schematic layout of the integrate-and-firecomponents. The graph (right) depicts voltage vs time for a neuronstimulated by a constant current.
	}\label{fig:2}	
\end{figure}

In the simulation comparing the EKF and SMA filters for estimating position and velocity, we use real eye-gaze tracking data and filtered estimates from the dataset provided in \cite{Judd2009} for our implementation. The main results of this comparison are presented in Fig. \ref{fig:2}.

For the real-world data, the RMSE values for the EKF and SMA are:
\begin{itemize}
	\item[$\bullet$]  RMSE (EKF) = 1.21 for position and 0.18 for velocity.  
	\item[$\bullet$] RMSE (SMA) = 66.3 for position and 62.3 for velocity.  
\end{itemize}
The EKF again outperforms the SMA in terms of accuracy. This is particularly evident in the position estimation, where the EKF has a significantly lower RMSE. Although the SMA is still effective at reducing noise, it introduces a more significant error, especially when the data shows more abrupt movements or rapid changes. Furthermore, the real-world data often contains more complex noise patterns, including sensor artifacts, calibration errors, and other real-world inaccuracies. The EKF’s recursive estimation method allows it to better handle these complexities by continuously adjusting its estimates based on incoming data, while the SMA’s fixed-window averaging method may struggle with more dynamic or irregular noise.

%Summary of Stochasticity**
%
%The stochastic components of the algorithm are present primarily in the **process noise** \( Q \) and **measurement noise** \( R \), which are used in the EKF to account for uncertainty in both the state evolution and the observations. These noises are modeled as random variables, and the EKF propagates these uncertainties over time, adjusting the state estimate accordingly.
%
%In summary:
%- **Process Noise** \( Q \) models uncertainty in how the system evolves.
%- **Measurement Noise** \( R \) models uncertainty in the observations.
%- The **EKF** integrates these stochastic elements to produce an estimate that minimizes the effect of noise on the state estimate.
%- **RMSE** is used to evaluate the accuracy of the estimates but is influenced by the level of stochastic noise in the system.

Thus, the EKF accounts for the system's stochastic nature using noise covariances, which allows the filter to estimate the state despite uncertainty in the process and measurements.

\section{Conclusions and discussion}
We have proposed an EKF to smooth gaze data collected during eye-tracking experiments. Our numerical results
have shown that the EKF effectively enhances the quality of the gaze signal by reducing high-frequency noise and compensating for data loss during blinks. Through simulations, we observe that the EKF smooths both position and velocity estimates, providing a more consistent and accurate representation of the subject's gaze trajectory. Quantitative metrics, such as root mean square error, show that the EKF outperforms traditional filtering methods, like the simple moving average, in terms of both noise reduction and accuracy.
The proposed approach represents a significant advancement in gaze tracking technology, offering both improved data quality and robustness to real-world challenges. This methodology can be extended to other relevant applications in human-computer interaction, neuroengineering, and behavioral research, where accurate real-time tracking of movement is crucial. Future work will explore further enhancements to the model, including applying deep learning and machine learning techniques \cite{Algabri2024deep}, and its potential integration with multimodal tracking systems to achieve even higher levels of precision and reliability.

	\bibliographystyle{splncs04}
\bibliography{mybibn-old}

\begin{thebibliography}{10}
\providecommand{\url}[1]{\texttt{#1}}
\providecommand{\urlprefix}{URL }
\providecommand{\doi}[1]{https://doi.org/#1}

\bibitem{Abaspourazad2024}
Abbaspourazad, H., Erturk, E., Pesaran, B., Shanechi, M.M.: Dynamical flexible
  inference of nonlinear latent factors and structures in neural population
  activity. Nature Biomedical Engineering  \textbf{8}(1),  85--108 (2024)

\bibitem{Algabri2024deep}
Algabri, R., Abdu, A., Lee, S.: Deep learning and machine learning techniques
  for head pose estimation: a survey. Artificial Intelligence Review
  \textbf{57}(10),  1--66 (2024)

\bibitem{Belykh2023}
Belykh, I., Kuske, R., Porfiri, M., Simpson, D.J.: Beyond the bristol book:
  Advances and perspectives in non-smooth dynamics and applications. Chaos: An
  Interdisciplinary Journal of Nonlinear Science  \textbf{33}(1) (2023)

\bibitem{Judd2009}
Judd, T., Ehinger, K., Durand, F., Torralba, A.: Learning to predict where
  humans look. In: 2009 IEEE 12th International Conference on Computer Vision.
  pp. 2106--2113 (2009)

\bibitem{Khodarahmi2023}
Khodarahmi, M., Maihami, V.: A review on kalman filter models. Archives of
  Computational Methods in Engineering  \textbf{30}(1),  727--747 (2023)

\bibitem{Law2015data}
Law, K., Stuart, A., Zygalakis, K.: Data assimilation. Cham, Switzerland:
  Springer  \textbf{214} (2015)

\bibitem{Meo2024}
Meo, M.M., Iaconis, F.R., Del~Punta, J.A., Delrieux, C.A., Gasaneo, G.:
  Multifractal information on reading eye tracking data. Physica A: Statistical
  Mechanics and its Applications  \textbf{638},  129625 (2024)

\bibitem{Niehorster2020}
Niehorster, D.C., Zemblys, R., Beelders, T., Holmqvist, K.: Characterizing gaze
  position signals and synthesizing noise during fixations in eye-tracking
  data. Behavior Research Methods  \textbf{52},  2515--2534 (2020)

\bibitem{Pyrhonen2023}
Pyrh{\"o}nen, L., Jaiswal, S., Mikkola, A.: Mass estimation of a simple
  hydraulic crane using discrete extended kalman filter and inverse dynamics
  for online identification. Nonlinear Dynamics  \textbf{111}(23),
  21487--21506 (2023)

\bibitem{Swari2021}
Swari, M.H.P., Handika, I.P.S., Satwika, I.K.S.: Comparison of simple moving
  average, single and modified single exponential smoothing. In: 2021 IEEE 7th
  Information Technology International Seminar (ITIS). pp.~1--5. IEEE (2021)

\bibitem{Toivanen2016}
Toivanen, M.: An advanced kalman filter for gaze tracking signal. Biomedical
  Signal Processing and Control  \textbf{25},  150--158 (2016)

\bibitem{Wu2002neural}
Wu, W., Black, M., Gao, Y., Serruya, M., Shaikhouni, A., Donoghue, J.,
  Bienenstock, E.: Neural decoding of cursor motion using a kalman filter.
  Advances in neural information processing systems  \textbf{15} (2002)

\bibitem{Yadav2023}
Yadav, S., Saha, S.K., Kar, R.: An application of the kalman filter for eeg/erp
  signal enhancement with the autoregressive realisation. Biomedical Signal
  Processing and Control  \textbf{86},  105213 (2023)

\bibitem{Yang2023outlier}
Yang, B., Huang, J.: Outlier-robust gaze signal filtering framework based on
  eye-movement modality recognition and set-membership approach. IEEE
  Transactions on Biomedical Engineering  \textbf{70}(8),  2463--2474 (2023)

\bibitem{Zulhakim2023}
Zulhakim, A.M., Abdullah, W.F.H., Halim, I.S.A., Mamat, R.B.H., Muslan, M.I.A.,
  Bakar, A.Z.A.: Smoothing sensor data in a controlled iot framework with
  moving averages. In: 2023 IEEE Regional Symposium on Micro and
  Nanoelectronics (RSM). pp. 86--89. IEEE (2023)

\end{thebibliography}

\end{document}